\documentclass[final,3p,times,12pt]{elsarticle}

\usepackage{lineno}
\usepackage{amsmath}
\usepackage{amssymb}
\usepackage{amsthm}
\usepackage{cases}
\usepackage{url}
\usepackage{parskip}
\usepackage[pdfpagelabels,plainpages=false]{hyperref}
\hypersetup{ linkcolor =blue,
 citecolor = cyan,
 urlcolor = blue,
colorlinks = true, bookmarksopen = true, bookmarksnumbered = true,
bookmarksopenlevel = {1}, linktocpage = true,
pdfstartview = FitH,
pdftitle = {Mixed Recurrences, interlacing properties and bounds  of zeros},
pdfsubject
= {article},
pdfauthor = {JJ},
pdfcreator = {pdftextify}, pdfproducer = {pdftextify}}

\newtheorem{theorem}{Theorem}
\newtheorem{corollary}[theorem]{Corollary}
\newtheorem{lemma}[theorem]{Lemma}
\newtheorem{proposition}[theorem]{Proposition}

\newcommand{\nn}{\mathbb{N}}

\newtheorem{remark}[theorem]{Remark}
\def\a{\alpha}
\def\b{\b}
\newcommand{\rr}{\mathbb{R}}

\def\N{\mathbb{N}}

\def\b{\beta}

\newcommand{\comment}[1]{}
\journal{ }
\begin{document}
	
\begin{frontmatter}
\title{On zeros of quasi-orthogonal Meixner polynomials}

\author[label1]{A.S.~Jooste}\ead{alta.jooste@up.ac.za}
\author[label2]{K.~Jordaan\corref{cor1}}\ead{jordakh@unisa.ac.za}
\cortext[cor1]{Corresponding author}

\address[label1]{Department of Mathematics and Applied Mathematics, University of Pretoria, Pretoria 0002, South Africa}
\address[label2]{Department of Decision Sciences, University of South Africa, Pretoria 0002, South Africa}

\begin{abstract}
For each fixed value of $\beta$ in the range $-2<\b<-1$ and $0<c<1$, we investigate interlacing properties of the zeros of polynomials of consecutive degree for $M_{n}(x;\b,c)$ and $M_k(x,\b+t,c)$, $k\in\{n-1,n,n+1\}$ and $t\in\{0,1,2\}$. We prove the conjecture in \cite{DriverJooste} on a lower bound for the first positive zero of the quasi-orthogonal order $1$ polynomial $M_n(x;\b+1,c)$ and identify upper and lower bounds for the first few zeros of quasi-orthogonal order $2$ Meixner polynomials $M_n(x;\b,c)$. We show that a sequence of Meixner polynomials $\{M_n(x;\b,c)\}_{n=3}^{\infty}$ with $-2<\b<-1$ and $0<c<1$ cannot be orthogonal with respect to any positive measure by proving that the zeros of $M_{n-1}(x;\b,c)$ and $M_{n}(x;\b,c)$ do not interlace for any $n\in\mathbb{N}_{\geqq 3}.$
\end{abstract}


\end{frontmatter}

\section{Introduction }
Consider a sequence $\{p_n\}_{n=0}^{\infty}$ of monic orthogonal polynomials with zeros $x_{1,n}<x_{2,n}<\dots<x_{n,n},$ satisfying a  three-term recurrence relation
\begin{equation}p_{n}(x)=\label{TTRR1}(x-B_n)p_{n-1}(x)-C_np_{n-2}(x),\end{equation}  where $B_n$ and $C_n>0$ do not depend on $x$, $p_{-1}\equiv0$ and $p_0\equiv1$. A consequence of orthogonality is that each polynomial $p_n(x)$ has $n$ real, distinct zeros in  $(a,b)$ and the zeros of $p_n$ and $p_{n-1}$ interlace as follows \begin{equation*}x_{1,n}<x_{1,n-1}<x_{2,n}<\dots <x_{n-1,n}<x_{n-1,n-1}<x_{n,n}.\end{equation*}
 Further, the interlacing of zeros of $p_n$ and $q_m, m\le n-2,$ referred to as Stieltjes interlacing, means that there exist $m$ open intervals with endpoints at successive zeros of $p_n$, each of which contains exactly one zero of $q_m.$

 In order for orthogonality conditions to hold, we often need restrictions on the parameters of the classical orthogonal polynomials and when the parameters deviate from these restricted values in an orderly way, the zeros may depart from the interval of orthogonality in a predictable way. This phenomenon can be explained in terms of the concept of quasi-orthogonality. The sequence of polynomials  $\{Q_{n}\}_{n=r+1}^N,$ where each polynomial $Q_{n}$ is of exact degree $n$ and $N\in\mathbb{N}\cup \{\infty\},$ is discrete quasi-orthogonal of order $r\in\mathbb{N}$, where $0\leq r<n,$ with respect to the weight function $w(x)>0$ on $[a,b]$ if (cf. \cite{chiharaboek})
 \begin{equation}\label{quasi_orth}\sum_{i=0}^{N-1} (x_i)^m q_n(x_i)w(x_i)\begin{cases}=0,~\mbox{for} ~m\in\{0,1,\dots,n-r-1\},\\\neq0,~\mbox{for} ~m=n-r.\end{cases}\end{equation} Note that quasi-orthogonal polynomials $Q_n$ are only defined for $n\in\{r+1,r+2,\dots\}.$ It is clear that when $r=0$ in \eqref{quasi_orth}, the sequence  $\{Q_{n}\}_{n\geq 0}$ is discrete orthogonal with respect to $w(x)$ on $[a,b]$.

Quasi-orthogonality was first studied by Riesz \cite{riesz}, followed by Fej\'{e}r \cite{fejer}, Shohat \cite{Shohat},  Chihara \cite{Chihara}, Dickinson \cite{Dickinson}, Draux \cite{Draux_1990}, Maroni \cite{maroni} and Joulak \cite{Joulak_2005}. The quasi-orthogonality of  Jacobi, Gegenbauer and Laguerre sequences was discussed in \cite{BDR}, the quasi-orthogonality of  Meixner sequences in \cite{JJT} and of Meixner-Pollaczek, Hahn, Dual-Hahn and Continuous Dual-Hahn sequences in \cite{JJJ}. Quasi-orthogonality of Pseudo-Jacobi polynomials was considered in \cite{JT} while \cite{TJK} deals with quasi-orthogonality of some hypergeometric and $q$-hypergeometric polynomials. 

Quasi-orthogonal polynomials are characterised by the following property:
\begin{lemma}[\cite{BDR, Chihara}]\label{Brez}
  Let $\{P_n\}_{n\geq 0}$ be a family of orthogonal polynomials on $[a,b]$ with respect to the weight function $w(x)>0$. A necessary and sufficient condition for a polynomial $Q_{n,k}$ of degree $n$ to be quasi-orthogonal of order $k\le n-1$ with respect to $w$  on $[a,b]$, is that
  \begin{equation*}
    Q_{n,k}(x)=\sum_{i=0}^k a_{n,i}P_{n-i}(x),~ a_{n,0}a_{n,k}\neq 0.
  \end{equation*}
\end{lemma}

It is known that for a quasi-orthogonal polynomial of order $k$ with $n$ distinct real roots, there exists a quadrature formula valid for all polynomials of degree $\leq$ $2n-k-1$ based on these zeros, provided that the weight in the quadrature formula does not vanish (see, for example, \cite{Freud, Michelli, Xu}). $n$-point quadrature rules with one or two fixed nodes achieve the highest degree of accuracy when the remaining nodes are the zeros of a polynomial of degree $n$ in a sequence of quasi-orthogonal polynomials (cf. \cite{Beckermann,bultheel,bultheel2}). Interlacing properties of zeros of polynomials of consecutive degree were used in \cite{Lubinsky} to provide substantially shorter proofs for generalisations of quadrature identities.  

Interlacing of the zeros was studied for various classical quasi-orthogonal polynomials, namely  quasi-orthogonal order $1$ Gegenbauer polynomials in \cite{ DriverMuldoonGeg}, Laguerre polynomials in \cite{DriverMuldoonLag}, Jacobi polynomials in \cite{ DriverJordaan},  Meixner polynomials in \cite{DriverJooste},  Pseudo-Jacobi polynomials in \cite{JT}, $q$-Laguerre in \cite{KG1} and little $q$-Jacobi in \cite{KG}, and quasi-orthogonal  order $2$ Jacobi polynomials in \cite{DriverJordaanJac2}.

In this paper, we investigate zeros of quasi-orthogonal Meixner polynomials with a particular emphasis on bounds for and interlacing properties of the zeros. Monic Meixner polynomials can be defined in terms of the $_2F_1$ hypergeometric series  by  (cf. \cite[(9.10.1)]{KLS})\begin{align}\label{MMdef}M_n(x;\b,c)&=\left(\frac{c}{c-1}\right)^n \left(\b \right)_n \, _2F_1\left(-n,-x;\b ;1-\frac{1}{c}\right)\\&=\nonumber \left(\frac{c}{c-1}\right)^n(\b)_n\phantom{ }~
\sum_{k=0}^{n}\frac{(-n)_k(-x)_k(1-\frac 1c)^k}{(\b)_kk!},\phantom{ }~ \b,c\in
\rr, ~\b\neq-1,-2,\dots,-n+1, ~c\neq 0\end{align}
where $(~ )_n$ is Pochhammer's symbol defined by
\begin{eqnarray*}(\a)_n& =&(\a)(\a+1)\cdots(\a+n-1)~\mbox{for}~
n\geq1\\
(\a)_0&=&1 ~\mbox{when}~ \a\neq0.\end{eqnarray*}
The three term recurrence relation satisfied by monic Meixner polynomials \cite[(9.10.4)]{KLS} is 
\begin{equation}\label{MM3trr}M_{n}(x;\b,c)=\left(x+\frac{\b  c+c n-c+n-1}{c-1}\right) M_{n-1}(x;\b,c)-\frac{c (n-1) (\b +n-2)}{(c-1)^2} M_{n-2}(x;\b,c)\end{equation}
with initial conditions $M_{-1}(x;\b,c)=0$ and $M_{0}(x;\b,c)=1$. Meixner polynomials for various parameter values, including those for Krawtchouk polynomials, are discussed in \cite{JJT}. 

 In Section \ref{Meixner}, we briefly recap and expand on some results for the zeros of Meixner polynomials when $\b>0$ and $c\in(0,1)$. In Section \ref{qo1}, we focus on the parameter ranges $-\b,c\in(0,1)$, for which the Meixner polynomials are discrete quasi-orthogonal of order $1$ and resolve a conjecture posed in \cite{DriverJooste}. We consider the case when $\b\in(-2,-1)$ and $c\in(0,1)$ for which Meixner polynomials are discrete quasi-orthogonal of order $2$ in Section \ref{qo2}.

\section{Meixner polynomials $M_n(x;\b,c)$, $\b>0$, $0<c<1$}\label{Meixner}

The weight function \begin{equation}\label{Mweight}w(x)=\frac{c^x(\b)_x}{x!}\end{equation} is positive when $\b,c>0$ and, by the ratio test, the moments associated with the weight \eqref{Mweight} exist when $c\in[0,1)$. 

The sequence $\{M_n(x;\b,c)\}_{n=0}^{\infty}$ satisfies the
discrete orthogonality relation (cf. \cite[(9.10.2)]{KLS})
\begin{equation*}\sum_{x=0}^{\infty}\frac{c^x(\b)_x}{x!}M_m(x;\b,c)M_n(x;\b,c)=\frac{c^{-n}n!}{(\b)_n(1 - c)^{\b} }\phantom{}~\delta_{mn}\end{equation*} when $\b>0$ and $c\in(0,1)$.
For $0<c<1$ and $\b>0$, the zeros are all distinct, positive and monotonically increasing as $\b$ increases (cf. \cite[Thm 7.1.2]{Ismail}). Further, there exists a point of the support of the measure between any two consecutive zeros of a discrete orthogonal polynomial. The minimum distance between consecutive zeros of Meixner polynomials is greater than $1$ \cite[Thm 1]{KZ}. 

The following lemma will be used to prove one of our results.
\begin{lemma}\label{lemrem}
Let $c\in(0,1)$, $\b>0$ and $n\in\nn$. Denote the zeros of $M_n(x;\b,c)$ by $0<y_{1,n}<y_{2,n}<\dots<y_{n,n}$, then
\begin{itemize}
 \item[(i)] $y_{1,1}=\frac{\b c}{1-c} \le \b$ when $c\in(0,0.5]$;
    \item[(ii)] $y_{1,n}<1<y_{2,n}$ when $n\ge\frac{\b c}{1-c}$.
\end{itemize}
    \end{lemma}
    \begin{proof} \begin{itemize}
        \item[(i)] A simple calculation, using either \eqref{MMdef} or \eqref{MM3trr}, shows that $M_1(x;\b,c)=x+\frac{\b c}{c-1}$ and hence $y_{1,1}=\frac{\b c}{1-c}\leq \b$ for $c\in(0,0.5].$ 
        \item[(ii)] For $0<c<1$ and $\b>0$, we have $\displaystyle \frac{M_n(1;\b,c)}{M_n(0;\b,c)}=\frac{\b  c+c n-n}{\b  c}\le 0$ for $n\ge\frac{\b c}{1-c}$. See \cite[Remark 2, p.127 and p.131]{KZ} for an alternate proof using a difference equation satisfied by the polynomial. \end{itemize} 
   \end{proof}

\section{Meixner polynomials $M_n(x;\b,c)$, $-1<\b<0$, $0<c<1$}\label{qo1}
It follows from the contiguous relation for hypergeometric functions \cite[eqn. (20)]{JJT} that monic Meixner polynomials satisfy \begin{equation}\label{o1rel} M_n(x;\b,c)=M_n(x;\b+1,c)-\frac{nc}{c-1}M_{n-1}(x;\b+1,c)\end{equation} and hence, for $-1<\b<0$, the sequence $\{M_n(x;\b,c)\}_{n=2}^{\infty}$ is quasi-orthogonal of order $r=1$ by Lemma \ref{Brez}. For $n\geq2$, $0<c<1$ and $-1<\b<0,$ the zeros of $M_n(x;\b,c)$ are all distinct, the smallest zero is  always negative and the remaining $(n-1)$ zeros are all positive (cf. \cite[Thm 2.1(i)]{DriverJooste}).

The interlacing of zeros of polynomials within the sequences of quasi-orthogonal order $1$ Meixner polynomials $\{M_n(x;\b,c)\}^{\infty}_{n=2}$ characterised by $-\b,c\in(0,1)$ is discussed in \cite{DriverJooste} as well as the interlacing of zeros of quasi-orthogonal Meixner polynomials $M_n(x;\b,c)$ with the zeros of their nearest orthogonal counterparts $M_l(x;\b+k,c)$, $l,n\in\nn$, $k\in\{1,2\}.$ 
Driver and Jooste \cite{DriverJooste} conjectured that, if $z_{1,n}<0<z_{2,n}<\dots<z_{n,n}$ are the zeros of quasi-orthogonal order $1$ $M_n(x;\b,c)$, then $z_{2,n}>1$.  

\begin{proposition}\label{prop}
Suppose that $\{M_n(x;\b,c)\}_{n\geq 2}$ is a sequence of quasi-orthogonal order $1$ Meixner polynomials with $-\b,c\in(0,1)$. Denote the zeros of $M_n(x;\b,c)$ by $\{z_{i,n}\}_{i=1}^{n}$ in ascending order, then
    $$z_{1,n}<0<1<z_{2,n}.$$
   \end{proposition}
\begin{proof} Let $-\b,c\in(0,1)$, then $\displaystyle \frac{M_n(1;\b,c)}{M_n(0;\b,c)}=\frac{\b  c+c n-n}{\b  c}>0$ for $n>\frac{\b c}{1-c}$. Since $\frac{\b c}{1-c}<0$ when   $-\b,c\in(0,1)$, the polynomial $M_n(x;\b,c)$ has either no zeros or an even number of zeros in $(0,1)$, for all $n\ge 0$. From \cite[Thm 2.1(i)]{DriverJooste},  $z_{1,n}<0$ for all $n$, and by \cite[Remark (2), p.482] {DriverJooste}, $z_{3,n}>1$. It therefore is impossible for an even number of zeros of $M_n(x;\b,c)$ to lie in $(0,1)$ and we deduce that there are no zeros of $M_n(x;\b,c)$ in $(0,1)$. Therefore $z_{1,n}<0<1<z_{2,n},$ which proves  \cite[Conjecture I]{DriverJooste}.
 \end{proof}

The following results are refinements of results in \cite{DriverJooste} that relied on the conjecture as an assumption. We state them here without proof.

\begin{theorem}\label{DJTh2.2}(cf. \cite[Thm 2.2]{DriverJooste})
Fix $n,\b$ and $c$ where  $n\in\N_{\geqq 3}$ and $-\b,c\in(0,1)$, then the zeros of $xM_{n-1}(x;\beta, c)$ and $M_{n}(x;\beta, c)$ interlace.
\end{theorem}
\begin{corollary} \label{DJCor2.3}(cf. \cite[Cor. 2.3]{DriverJooste})
Let $-\b,c\in(0,1)$ and $n\in\N_{\geqq 3}.$ Then, if $z_{1,n}<0<z_{2,n}<\dots<z_{n,n}$ denote the zeros of $M_{n}(x;\beta, c)$,
 the negative zero of $M_n(x;\beta,c)$ increases with $n.$  Moreover, $z_{1,1}=\frac{\b c}{1-c}$ is a lower bound for the negative zero of $M_n(x;\b,c)$ for each $n\in\N.$
\end{corollary}
\begin{theorem}\label{DJTh2.4} (cf. \cite[Thm 2.4]{DriverJooste})
Fix $n,\b$ and $c$ where $n\in\N_{\geqq 4}$, $-\b,c\in(0,1)$. Assume that $M_n(x;\beta,c)$ and $M_{n-2}(x;\beta,c)$ do not have common zeros, i.e., $M_n\left(\frac{\beta c+(c+1)(n-1)}{1-c};\beta,c\right)\neq 0$ (or, equivalently, $M_{n-2}\left(\frac{\beta c+(c+1)(n-1)}{1-c};\beta,c\right)\neq 0$).  Then the positive zeros of $\left(x-\frac{\beta c+(c+1)(n-1)}{1-c}\right)M_{n-2}(x;\beta, c),$ interlace with the positive zeros of $M_n(x;\beta,c).$
\end{theorem}
\begin{theorem}\label{DJTh3.1} (cf. \cite[Thm 3.1]{DriverJooste})
Fix $n,\b$ and $c$ where $n\in\N$ and $-\b,c\in(0,1)$. Assume  that $M_n(x;\beta,c)$ and $M_{n-2}(x;\beta+1,c)$ do not have common zeros i.e., $M_n\left(\frac{\b c + n-1}{1-c};\beta,c\right)\neq 0$ (or, equivalently, $M_{n-2}\left(\frac{\b c + n-1}{1-c};\beta,c\right)\neq 0$). Then, for each $n\in\N_{\geqq 3},$ the zeros of $\left(x-\frac{\b c + n-1}{1-c}\right)M_{n-2}(x;\beta+1, c)$  interlace with the $(n-1)$ positive zeros of $M_n(x;\beta,c).$
\end{theorem}

The final result in this section considers a case not addressed in \cite{DriverJooste}, namely when both the degree and $\b$ are increasing. 
\begin{theorem}\label{Eq55} Suppose that $\{M_n(x;\b,c)\}_{n\geq 2}$ is a sequence of quasi-orthogonal order $1$ Meixner polynomials with $-\b,c\in(0,1)$. Let $C_n=-\b+\frac{(n+1)c}{1-c}$ and assume that $M_n(C_n;\b,c)\neq 0$, then the $(n+1)$ zeros of $(x-C_n)M_{n}(x;\b, c)$ interlace with the $(n+1)$ zeros of $M_{n+1}(x;\b+1, c)$.
\end{theorem}
 
\begin{proof}

Let $-\b,c\in(0,1)$ be fixed. We note that $C_n>-\b$ for these values of $\b$ and $c$, furthermore $z_{1,n}<0<1<z_{2,n}$ for all values of $n$. Consider the equation
\begin{equation}\label{555} \frac{  \b+n}{1 - c}M_{n}(x;\b,c) = \left(x -C_n\right) M_{n}(x;\b+1,c)-M_{n+1}(x;\b+1,c),\end{equation} which can be verified by comparing coefficients of $x^n$. Since $M_n(x;\b+1,c)$ and $M_{n+1}(x;\b+1,c)$ are polynomials of consecutive degree in an orthogonal sequence, their zeros are interlacing and they have no common zeros. It follows from \eqref{555} that $M_n(x;\b,c)$ and $M_{n+1}(x;\b+1,c)$ can have at most one common zero, namely at $x=C_n$. If they do not have any common zeros, i.e. $M_n(C_n;\b,c)\neq 0$ (equivalently, $M_{n+1}(C_n;\b+1,c)\neq 0$), evaluating (\ref{555}) at $y_{k,n+1}$ and  $y_{k+1,n+1}$, $k\in\{1,2,\dots,n\}$, two consecutive zeros of  $M_{n+1}(x;\b+1, c)$, we obtain
\begin{equation}\label{2}
\frac{\left(\frac{\b+n}{ 1-c}\right)^2 }{M_n(y_{k,n+1};\b +1,c)M_n(y_{k+1,n+1};\b+1 ,c)}= \frac{\left(y_{k,n+1}-C_n\right)\left(y_{k+1,n+1}-C_n\right)}{ M_{n}(y_{k,n+1};\b ,c) M_{n}(y_{k+1,n+1};\b ,c)}.\end{equation}
The zeros of the polynomials in an orthogonal sequence interlace, therefore $M_n(x;\b+1,c)$ differs in sign at the zeros of $M_{n+1}(x;\b+1, c)$ and $M_n(y_{k,n+1};\b +1,c)M_n(y_{k+1,n+1};\b+1 ,c)<0$. The left hand side of (\ref{2}) is therefore negative. 

Suppose $C_n \in (y_{k,n+1},y_{k+1,n+1})$ for some $k\in\{1,2,\dots,n\}.$ Then in this single interval containing $C_n$, we will have  $\left(y_{k,n+1}-C_n\right)\left(y_{k+1,n+1}-C_n\right)<0$ and for the right hand side of (\ref{2}) to be negative, $M_{n}(y_{k,n+1};\b ,c) M_{n}(y_{k+1,n+1};\b ,c)>0$ and we  deduce that in this interval containing $C_n$, there will be no zeros of $ M_{n}(x;\b ,c).$ In each one of the other (n-1) intervals $(y_{k,n+1},y_{k+1,n+1}),k\in\{1,2,\dots,n\},$ not containing $C_n,$     there will be a zero of $ M_{n}(x;\b ,c).$ 

The point $C_n$ cannot lie to the left of $y_{1,n+1},$ since, for $-\b,c\in(0,1)$, we have  $y_{1,1}=\frac{(\b+1)c}{1-c}$ (see Lemma \ref{lemrem}) and  $$ \frac{C_n}{y_{1,1}}=\frac{\frac{c (n+1)}{1-c}-\b}{\frac{(\b+1)c}{1-c}}>1,$$
for $n>\frac{\b}{c},$ but $\frac{\b}{c}<0,$ therefore $C_n>y_{1,1}>y_{1,n}>y_{1,n+1}$ for all $n\ge 1$, which means $C_n \notin (0,y_{1,n+1})$.

Suppose $C_n>y_{n+1,n+1}$. Then $\left(y_{k,n+1}-C_n\right)\left(y_{k+1,n+1}-C_n\right)>0$ for each $k\in\{1,2,\dots,n\}$ and for the right hand side of (\ref{2}) to be negative, we need 
$M_{n}(y_{k,n+1};\b ,c) M_{n}(y_{k+1,n+1};\b ,c)<0$ for each $k\in\{1,2,\dots,n\}$, and this can only be true if there is a zero of $M_n(x;\b,c)$ in each one of the $n$ intervals with endpoints at the $(n+1)$ zeros of $M_{n+1}(x;\b+1,c).$ This leads to a contradiction, since 
 $ M_{n}(x;\b ,c)$ has $n$ zeros of which the first one is negative, which means only $(n-1)$ zeros of $ M_{n}(x;\b ,c)$ are available to fill the $n$ gaps between the $(n+1)$ positive zeros of $M_{n+1}(x;\b+1, c)$. We thus have $C_n < y_{n+1,n+1}$ and $C_n$ fills the $n$th gap. 
\end{proof}

\section{Meixner polynomials $M_n(x;\b,c)$, $-2<\b<-1$, $0<c<1$}\label{qo2}
 Iterating \eqref{o1rel}, we see that monic Meixner polynomials satisfy \cite[Thm 7]{JJT} 
$$M_{n}(x;\b ,c)=M_{n}(x;\b+2 ,c)+2n\frac{c }{1-c}M_{n-1}(x;\b+2 ,c)+n(n-1)\left(\frac{c}{c-1}\right)^2 M_{n-2}(x;\b+2 ,c).$$ This implies that, for $\b \in(-2,-1)$, the right hand side is a linear combination of three terms in a sequence of orthogonal polynomials and it follows from Lemma \ref{Brez} that the sequence $\{M_n(x;\b,c)\}_{n=3}^{\infty}$ is quasi-orthogonal of order $r=2$  when $c\in (0,1)$ and $\b \in (-2,-1)$. 

In what follows, we assume that $-2<\b<-1$, $0<c<1$ and use the following notation: For each $n\in\nn$ and $-2<\b<-1$ and $0<c<1$ fixed,  we will indicate the zeros of quasi-orthogonal order $2$ polynomials $M_{n}(x;\b, c)$ by $x_{1,n}<x_{2,n}<\dots<x_{n,n}$, the zeros of quasi-orthogonal order $1$ polynomials $M_{n}(x;\b+1, c)$ by $z_{1,n}<z_{2,n}<\dots<z_{n,n}$ and the zeros of orthogonal polynomials $M_{n}(x;\b+2, c)$ by $y_{1,n}<y_{2,n}<\dots<y_{n,n}$. By orthogonality, we have $$y_{1,n}<y_{1,n-1}<y_{1,n-2}<\cdots<y_{1,1}=\left(\frac{c}{1-c}\right) (\b+2)$$
and, if $c\in \left(0,0.5\right]$, it follows from Lemma \ref{lemrem}(i) that  
  $$y_{1,n}<y_{1,n-1}<\cdots<y_{1,1}\le \b+2 <1.$$

\begin{lemma}\label{00} Suppose $\{M_{n}(x;\b,c)\}_{n=3}^\infty$ is a sequence of quasi-orthogonal order 2 Meixner polynomials with $\b \in (-2,-1)$ and $c\in(0,1)$. Let $\{x_{i,n}\}_{i=1}^{n}$ denote the zeros of $M_{n}(x;\b,c)$ and $\{y_{i,n}\}_{i=1}^{n}$ the zeros of $M_{n}(x;\b+2, c)$ in increasing order. If $n>\frac{\b}{c-1}$, the zeros of $M_{n}(x;\b,c)$ are all distinct, positive and interlace with the $(n-1)$ zeros of $M_{n-1}(x;\b +2,c)$ as follows:
\[x_{1,n}<y_{1,n-1}<x_{2,n}<\dots<y_{n-1,n-1}<x_{n,n}.\]
\end{lemma}
\begin{proof}
The fact that the zeros of $M_n(x;\b,c)$ are simple and positive for $n>\frac{\b}{c-1}$ was proved in \cite[Thm 7]{JJT} using \cite[Thms 8 and 9]{Joulak_2005} but also follows from earlier results due to Xu (cf. \cite[Thms 5.1 and 5.3]{Xu}). The interlacing of the zeros of $M_n(x;\b,c)$ and $M_{n-1}(x;\b+2,c)$ is a straightforward consequence of \cite[Thm 10]{Joulak_2005}.
\end{proof}

\begin{remark}
Lemma \ref{00} applies to $n\geq3$ and the case when $n=1$ was discussed in Lemma \ref{lemrem}. For $n=2$, we have $$M_2(x;\b,c)=x^2+\frac{1+c+2 \b  c}{c-1}x+\frac{\b  (\b +1) c^2}{(c-1)^2}$$ with zeros
$x=\displaystyle \frac{2 \beta  c+c+1\pm\sqrt{4 \beta  c+(c+1)^2}}{2(1-c)}.$
\begin{itemize}
\item [(i)] $M_2(x;\b,c)$ will have a double root, $x=\frac{2 \b c+c+1}{2 (1-c)},$ when $\b=-\frac{(c+1)^2}{4c}$ and if $\b\in(-2,-1)$, this implies that $c\in \left(3-2\sqrt2,1\right)$.
\item [(ii)]$M_2(x;\b,c)$ will have two distinct roots if $\b>-\frac{(c+1)^2}{4c}$. Since $-\frac{(c+1)^2}{4c}<-1$ for all $c\in(0,1)$, there are two cases to consider. When $-2<-\frac{(c+1)^2}{4c}$, we will have 2 distinct roots  when $c\in\left(3-2 \sqrt{2},1\right)$ and $-2<-\frac{(c+1)^2}{4c}<\b<-1$ and in the case when  $-\frac{(c+1)^2}{4c}<-2$, we will have two distinct real roots when $c\in\left(0,3-2 \sqrt{2}\right)$ and $\b\in(-2,-1).$
\item [(iii)] $M_2(x;\b,c)$ will have two pure imaginary roots when $\b<-\frac{(c+1)^2}{4c}$, i.e., when $c\in\left(3-2 \sqrt{2},1\right)$ and $-2<\b<-\frac{(c+1)^2}{4c}<-1$.
\end{itemize}
\end{remark}

Next, we provide upper and lower bounds for the first three zeros of a quasi-orthogonal order $2$ Meixner polynomial which we need to prove our main results. 

\begin{lemma} \label{Th00}
Suppose $\{M_{n}(x;\b,c)\}_{n=3}^\infty$ is a sequence of quasi-orthogonal order 2 Meixner polynomials with $\b \in (-2,-1)$ and $c\in(0,1)$. If $n >\frac{\b}{c-1}$ and $\{x_{i,n}\}_{i=1}^{n}$ denotes the zeros of $M_{n}(x;\b, c),$ then $$0<x_{1,n}<-\b-1<x_{2,n}<1<-\b<2<x_{3,n}.$$
\end{lemma}

 \begin{proof}
Fix $\b\in(-2,-1)$ and $c\in(0,1)$ and suppose $n \in\nn_{\geqq 3}$ with $n >\frac{\b}{c-1}$. Then $-\b$ and $-\b-1$ are not zeros of $M_{n}(x;\b, c)$, since  $ M_n(-\b-1;\b,c)=\left(\frac{1}{c-1}\right)^n (\b-c n+n) \frac{\Gamma (\b+n)}{\Gamma (\b+1)}=0$ only 
when $n=\frac{\b}{c-1}$ and  $M_n(-\b;\b,c)=\left(\frac{1}{c-1}\right)^n (\b)_n\neq 0$.
 
Let $y_{1,n}<y_{2,n}<\dots<y_{n,n}$ denote the zeros of $M_{n}(x;\b +2, c).$ From Lemma \ref{00}, we know that $x_{1,n}<y_{1,n-1}$. From Lemma \ref{lemrem}(ii), we know that $y_{1,n}<1$ for all $n>\frac{(\b+2)c}{1-c}$, which implies  $y_{1,n-1}<1$ when $n>\frac{(\b+2)c}{1-c}+1=\frac{\b c+c+1}{1-c}$. Since
$$\frac{\frac{\b c+c+1}{1-c}}{\frac{\b}{c-1}}=-\frac{\b  c+c+1}{\b }<1$$
for all $\b\in(-2,-1)$ and $c>-1$, i.e., also when $c\in(0,1)$, we have  $\frac{\b c+c+1}{1-c}<\frac{\b}{c-1}$ for all $\b\in(-2,-1)$ and $c\in(0,1)$. Consequently $x_{1,n}<y_{1,n-1}<1$ when $n>\frac {\b}{c-1}.$ Furthermore, since 
$$ \frac{M_n(1;\b,c)}{M_n(0;\b,c)}=\frac{\b c+c n-n}{\b  c}>0$$ when $n>\frac{\b c}{1-c}$,  the polynomial $M_n(x;\b,c)$, has either no zero or an even number of zeros in $(0,1)$ and this holds true for all $n\ge0 $, since $\frac{\b c}{1-c}<0$ if $0<c<1$ and $-2<\b<-1$. Therefore, in the case when all the zeros of  $M_{n}(x;\b, c),$ are real, i.e., when $n >\frac{\b}{c-1},$ we have 
 $0<x_{1,n}<x_{2,n}<1$.

 Now, for $0<c<1$, $-2<\b<-1$ and $n>\frac{\b}{c-1}$, consider
$$\displaystyle \frac{M_n(-\b-1;\b,c)}{M_n(0;\b,c)}=\frac{n+\b -c n}{\b c^n}<0$$  and $$\displaystyle \frac{M_n(-\b-1;\b,c)}{M_n(-\b;\b,c)}=\frac{n+\b-c n}{\b}<0.$$

Hence the polynomial $M_n(x;\b,c)$ has an odd number of zeros in each of the intervals $(0,-\b-1)$ and $(-\b-1,-\b)$, therefore we have 
$$0<x_{1,n}<-\b-1<x_{2,n}<1<-\b<2.$$
From the interlacing proved in Lemma \ref{00}, we know that, for all  $n>\frac{\b}{c-1},$ 
\begin{equation}\label{config}x_{2,n}<y_{2,n-1}<x_{3,n}<y_{3,n-1}<x_{4,n}.\end{equation} Since the zeros of the Meixner polynomials, in this case $y_{2,n-1}$ and $y_{3,n-1},$ are more than one unit apart, it follows naturally that there cannot be $3$  zeros of $M_n(x;\b,c)$ in the interval $(-\b-1,-\b)$, therefore $x_{3,n}>-\b$ and we have the configuration
 $$0<x_{1,n}<-\b-1<x_{2,n}<1<-\b<x_{3,n}.$$
Finally, we prove that $x_{3,n}>2$. From (\ref{config})
 and the fact that the zeros of Meixner polynomials are more than one unit apart, it follows that $x_{4,n}>y_{3,n-1}>2$ for $n>\frac{\b}{c-1}.$

 Consider \begin{align*}\frac{M_n(2;\b,c)}{M_n(-\b;\b,c)}&=c^{n-2}~\frac{ \beta  (\beta +1) c^2+n(c-1) \left[n(c-1)+2 \beta  c+c+1\right]}{\beta  (\beta +1)}.\end{align*}

The denominator as well as the first term in the numerator are positive for $-2<\b<-1$. Since $n>\frac{\b}{c-1}$ and $0<c<1$, we have that $n(c-1)<\b<-1$ and
\begin{equation*}
n(c-1)+2 \b c+c+1< (2\b+1)c< -c.\end{equation*} Hence the expression 
$n(c-1)(n(c-1)+2 \b c+c+1)$ is positive for $n>\frac{\b}{c-1},$ $0<c<1$ and $-2<\b<-1$
and it follows that $\displaystyle \frac{M_n(2;\b,c)}{M_n(-\b;\b,c)}>0$. This implies that there are either no zeros or an even number of zeros of $M_n(x;\b,c)$ in the interval $(-\b,2)$. The inequalities $x_{2,n}<1$, $x_{3,n}>-\b$ and $x_{4,n}>2$ imply that there cannot be any zeros of $M_n(x;\b,c)$ in $(-\b,2)$ and hence $x_{3,n}>2$.
 \end{proof}

The result in Lemma \ref{Th00} is illustrated in Table \ref{position}.

\begin{table}[ht] \centering \caption{The values of $x_{1,10}$, $x_{2,10}$ and $x_{3,10}$ for different values of $\b$ and $c$, in each case $n>\frac{\b}{c-1}$.}\label{position}

\begin{tabular}{c|c|c|c|c|c}
$-\b$&$c$&$\frac{\b}{c-1}$&$x_{1,10}$&$x_{2,10}$&$x_{3,10}$\\
\hline
$1.99$&$0.1$&$2.21$&$ 3.549 * 10^{-14}$&$0.999999999987$&$2.0000000038$\\
$1.99$&$0.5$&$3.98$&$1.351*10^{-6}$ &$0.999961$&$2.00059$\\
$1.99$&$0.8$&$9.5$& $0.001545$ & $0.990203$ &$2.02091$\\
$1.5$&$0.1$&$1.67$& $1.953*10^{-12}$ &$0.99999999904$&$2.00000054$\\
$1.5$&$0.5$&$3$& $0.0000535$ & $0.9982$& $2.0474$\\
$1.5$&$0.8$&$7.5$& $0.0405146$ & $0.811339$ &$2.93765$ \\
$1.01$&$0.1$&$1.12$& $1.230*10^{-13}$ & $0.999999999911$& $2.0000024$\\
$1.01$&$0.5$&$2.02$& $2.418*10^{-6}$ &$0.999893$& $2.1284$\\
$1.01$&$0.8$&$5.05$& $0.000973697$ & $0.994013$&$3.8303$\\
\end{tabular}\end{table}

\begin{theorem}\label{0} Suppose $\{M_{n}(x;\b,c)\}_{n=3}^\infty$ is a sequence of quasi-orthogonal order 2 Meixner polynomials with $\b \in (-2,-1)$ and $c\in(0,1)$. Let $\{x_{i,n}\}_{i=1}^{n}$ be the zeros of $M_{n}(x;\b, c)$, $\{y_{i,n}\}_{i=1}^{n}$ the zeros of $M_{n}(x;\b+2, c)$ in ascending order and $A_n=\frac{\b}{c-1}-\b-n-1$. 

\begin{itemize}
\item[(a)] If $n\ge \frac{\b}{c-1}-(\b+1)$, then $A_n\le 0$ and the zeros of $M_{n}(x;\b, c)$ and $M_{n}(x;\b+2, c)$ are interlacing:
\begin{eqnarray}\label{55}
x_{1,n}<y_{1,n}<x_{2,n}<y_{2,n}<\dots<x_{n,n}<y_{n,n}.\end{eqnarray}

\item [(b)]
There is at most one integer $n^*$ such that $\frac{\b}{c-1}-(\b+2)<n^*< \frac{\b}{c-1}-(\b +1)$ and  $0<A_{n^*}<1$.  When $n^* \in \left(\frac{\b}{c-1}, \frac{\b}{c-1}-(\b+1)\right)$, the zeros of $M_{n^*}(x;\b, c)$ and $M_{n^*}(x;\b+2, c)$ interlace as in (\ref{55}) if and only if $A_{n^*}<y_{1,n^*}$. 

\end{itemize}

\end{theorem}
\begin{proof}
Let $c\in(0,1)$ and $\b\in(-2,-1)$ be fixed.  Since $M_{n-1}(x,\b+2,c)$ and $M_n(x;\b+2,c)$ are polynomials of consecutive degree in an orthogonal sequences, their zeros interlace and they do not have common zeros. Consider
\begin{equation}\label{1}\left(\frac{\b }{ n}+1\right)M_n(x;\b ,c)=\left(\frac{\b}{n}+1-c\right)M_n(x;\b+2,c)+c \left(x-A_n\right) M_{n-1}(x;\b +2,c),\end{equation} which can be verified by comparing coefficients of $x^n$. It follows from \eqref{1} that $M_n(x;\b,c)$ and $M_n(x;\b+2,c)$ can have at most one common zero at $x=A_n$. Evaluating (\ref{1}) at $y_{k,n}$ and  $y_{k+1,n},k\in\{1,2,\dots,n-1\}$, consecutive zeros of  $M_{n}(x;\b+2, c)$, yields
\begin{equation}\label{3}
\frac{\left(\frac{\b }{ n}+1\right)^2\large{/}c^2}{  M_{n-1}(y_{k,n};\b +2,c) M_{n-1}(y_{k+1,n};\b +2,c)}=\frac{\left(y_{k,n}-A_n\right)\left(y_{k+1,n}-A_n\right)}{ M_n(y_{k,n};\b ,c)M_n(y_{k+1,n};\b ,c)}.\end{equation}

The polynomial $M_{n}(x;\b +2, c)$ belongs to an orthogonal sequence and since the zeros of two polynomials of consecutive degree, belonging to the same orthogonal sequence, interlace, $M_{n-1}(x;\b +2, c)$ will differ in sign at the two consecutive zeros $y_{k,n}$ and  $y_{k+1,n}$ of $M_{n}(x;\b +2, c)$ for each  $k\in\{1,2,\dots,n-1\}$, i.e.,
$M_{n-1}(y_{k,n};\b +2,c) M_{n-1}(y_{k+1,n};\b +2,c)<0$ and the left hand side of (\ref{3}) is negative.

\begin{itemize}
\item[(a)] Suppose $n\geq \frac{\b}{c-1}-(\b+1)$. It is clear that $A_n\le 0$ (and vice versa).  Furthermore, since $\b \in(-2,-1)$,  we also have  $n\geq \frac{\b}{c-1}-(\b+1)>\frac{\b}{c-1}$ and it follows from Lemma \ref{00} that the zeros of $M_{n}(x;\b, c)$ are positive and simple. Since the zeros of $M_n(x;\b,c)$ and $M_n(x;\b+2,c)$ are all positive and $A_n$, the only possible common zero, is negative, the two polynomials cannot have any common zeros in this case. Furthermore,
$A_n\notin(y_{k,n},y_{k+1,n})$ for each $k\in\{1,2,\dots,n-1\}$ which implies that the numerator on the right hand side of (\ref{3}) is positive. Since the  quotient on the right hand side of (\ref{3}) must be negative, $M_n(y_{k,n};\b ,c)M_n(y_{k+1,n};\b ,c)<0$ for each  $k\in\{1,2,\dots,n-1\}$, i.e., the polynomial $M_{n}(x;\b, c)$ differs in sign at the zeros of  $M_{n}(x;\b+2, c)$, therefore the  zeros of  $M_{n}(x;\b, c)$ and $M_{n}(x;\b+2, c)$ interlace. Since $n\geq \frac{\b}{c-1}-(\b+1)>\frac{(\b+2)c}{1-c}$ for $c\in(0,1)$ and $\b\in(-2,-1)$, the condition of Lemma \ref{lemrem}(ii) is met and we have $y_{2,n}>1$. Using Lemma \ref{Th00}, we deduce that the zeros interlace as follows:
\begin{eqnarray*}
0<x_{1,n}<y_{1,n}<x_{2,n}<1<y_{2,n}<\dots<x_{n,n}<y_{n,n}.\end{eqnarray*}
\item[(b)]
$\frac{\b}{c-1}-(\b+2)<n< \frac{\b}{c-1}-(\b +1)$ if and only if $0<A_n<1$ and since $$\frac{\b}{c-1}-(\b +1)-\left(\frac{\b}{c-1}-(\b +2)\right)=1,$$
$n$ lies in an interval of length $1$ and there is at most one value of $n$, say $n^*$, such that $0<A_{n^*}<1$. When $n^* \in \left(\frac{\b}{c-1}, \frac{\b}{c-1}-(\b+1)\right)$, all the zeros of $M_n(x;\b,c)$ are positive and simple. If $M_{n^*}(x;\b,c)$ and $M_{n^*}(x;\b+2,c)$ have a common zero at $A_{n^*}$, i.e. $M_{n^*}(A_{n^*};\b,c)=M_{n^*}(A_{n^*};\b+2,c)=0$, the zeros of $M_{n^*}(x;\b,c)$ and $M_{n^*}(x;\b+2,c)$ clearly do not interlace. Suppose $A_{n^*}<y_{1,n^*}$, then the zeros interlace as in \eqref{55} by the same argument as in (a). On the other hand, suppose the zeros interlace as in \eqref{55}. Since $0<A_{n^*}<1,$ the only possibility is that $A_{n^*}<y_{1,n^*}$.
\end{itemize}\end{proof}

\begin{remark}
In the limited number of cases, where $n\le \frac{\b}{c-1}-(\b+2)$, the two smallest zeros of $M_{n}(x;\b, c)$ are complex and numerical examples indicate that the $(n-2)$ remaining real zeros interlace with the largest $(n-2)$ zeros of the orthogonal polynomial $M_{n}(x;\b +2, c)$. 
\end{remark}

\noindent Next, we prove that if $0<c<1$ and $-2<\b<-1$ are fixed, the $(n-1)$ zeros of  $M_{n-1}(x;\b,c)$ do not interlace with the $n$ zeros of $M_{n}(x;\b,c)$ for any value of  $n \in\nn_{\geqq 4}$. 
  
\begin{theorem} \label{Th:2.1}
Suppose $\{M_{n}(x;\b,c)\}_{n=3}^\infty$ is a sequence of quasi-orthogonal order 2 Meixner
polynomials with $\b \in (-2,-1)$ and $c\in(0,1)$. Assume that $\b,c$ and $n$ are such that $M_{n}(x;\b,c)$ and $M_{n-1}(x;\b,c)$ have no common zeros. When $n-1>\frac{\b}{c-1},$ if $\{x_{i,n}\}_{i=1}^{n}$ are the zeros of $M_{n}(x;\b,c)$ in increasing order, the $n$ zeros of $(x+\b)M_{n-1}(x;\b,c)$ interlace  with  the $(n+1)$ zeros of $(x+\b+1)M_{n}(x;\b,c)$  as follows:
\begin{align}\label{2.2}0<x_{1,n} < x_{1,n-1}&< -\b-1 < x_{2,n-1}< x_{2,n} <-\b< x_{3,n} <x_{3,n-1}<x_{4,n}<\cdots\\&< x_{n-1,n}< x_{n-1,n-1}<x_{n,n}. \nonumber\end{align}	
\end{theorem}
\begin{proof}
Let $n-1>\frac{\b}{c-1}$, $-2<\b<-1$ and $0<c<1$. Evaluating the mixed recurrence relation 
\begin{equation}\label{3.15}\left(x+\frac{\b  (c-2)}{c-1}+n\right) M_n(x;\b,c)=-\frac{(\b +n-1) (\b -c n+n)}{(c-1)^2}M_{n-1}(x;\b,c)+(\b +x)_2  M_{n-1}(x;\b+2,c),\end{equation} 
at zeros $x_{i,n}$ and $x_{i+1,n}, i\in\{1,2,\dots,n-1\}$ of $M_{n}(x;\b,c)$, we obtain
\begin{align}\label{3.15ev34}\frac{\left(\frac{(\b +n-1) (\b -c n+n)}{(c-1)^2}\right)^2}{ M_{n-1}(x_{i,n};\b+2,c)M_{n-1}(x_{i+1,n};\b+2,c)}=\frac{(x_{i,n}+\b)(x_{i,n}+\b +1)(x_{i+1,n}+\b)(x_{i+1,n}+\b+1)}{M_{n-1}(x_{i,n};\b,c)M_{n-1}(x_{i+1,n};\b,c)}.\end{align}

We know from Lemma \ref{00} that, for $n>\frac{\b}{c-1}$, the zeros of $M_n(x;\b,c)$ interlace with the $(n-1)$ zeros of $M_{n-1}(x;\b+2,c)$, therefore $M_{n-1}(x_{i,n};\b+2,c)M_{n-1}(x_{i+1,n};\b+2,c)<0$ for $i\in\{1,\dots,n-1\}$ and the left hand side of (\ref{3.15ev34}) is negative.

For each $i\in\{3,4,\dots,n\}$, since $x_{i,n}>-\b$ from Lemma \ref{Th00}, the numerator on the right hand side  of (\ref{3.15ev34}) is positive, consequently $M_{n-1}(x_{i,n};\b,c)M_{n-1}(x_{i+1,n};\b,c)<0$ and this implies that each of the intervals $(x_{i,n},x_{i+1,n}), i\in\{3,4,\dots,n-1\}$ contains an odd number of zeros of $M_{n-1}(x;\b,c)$. Lemma \ref{Th00} with $n$ replaced by $n-1$ implies that $x_{1,n-1}<x_{2,n-1}<-\b$ for $n-1>\frac{\b}{c-1}$ and we see that there is exactly one of the remaining $(n-3)$ zeros of $M_{n-1}(x;\b,c)$ in each of the $(n-3)$ intervals $(x_{i,n},x_{i+1,n}), i\in\{3,4,\dots,n-1\}$, i.e.\begin{equation}
\label{lz}0<x_{1,n-1}< -\b-1 < x_{2,n-1}< 1<-\b< x_{3,n} <x_{3,n-1}<x_{4,n}<\dots< x_{n,n-1}< x_{n-1,n-1}<x_{n,n}    \end{equation}

To determine the relative positioning of the two smallest zeros of $M_{n-1}(x;\b,c)$ with respect to the two smallest zeros of $M_n(x;\b,c)$, consider \eqref{3.15ev34} for $i\in\{1,2\}$. Applying Lemma \ref{Th00} to the zeros of $M_n(x;\b,c)$, we see that the numerator on the right hand side in \eqref{3.15ev34} is negative in each case. Hence, $M_{n-1}(x_{1,n};\b,c)M_{n-1}(x_{2,n};\b,c)>0$ and $M_{n-1}(x_{2,n};\b,c)M_{n-1}(x_{3,n};\b,c)>0$, which implies that either no zeros or an even number of zeros of $M_{n-1}(x;\b,c)$ lie in each of the intervals $(x_{1,n},x_{2,n})$ and $(x_{2,n},x_{3,n})$. Considering \eqref{lz} and keeping in mind that $x_{1,n}< -\b-1 <x_{2,n}<1<-\b<x_{3,n}$ by Lemma \ref{Th00}, we see that the only possibility is for the zeros to be arranged as described in \eqref{2.2}. 

\end{proof}

\begin{remark}
\begin{itemize}
\item[]
\item[(i)]
Since the zeros of two quasi-orthogonal order $2$ Meixner polynomials of consecutive degree do not separate each other, the sequence of Meixner polynomials $\{M_{n}(x;\b,c)\}_{n=3}^{\infty}$ with $-2<\b<-1$ and $c\in(0,1)$ is not orthogonal with respect to any positive measure.   
\item[(ii)] Although the zeros of orthogonal Meixner polynomials, for $c\in(0,1)$ and $\b>0$, are monotonically increasing as $\b$ increases (cf. \cite[Thm 7.1.2]{Ismail}), the zeros of quasi-orthogonal order 1 and 2 Meixner polynomials are, in general, not monotonically increasing as $\b$ increases. For example, for $-1<\b<0$, the smallest zero is always negative, while for $-2<\b<-1$ and $n>\frac{\b}{c-1}$, the smallest zero is positive. The zeros of $M_n(x;\b,c)$ for a fixed value of $n$ and $c\in(0,1)$ and two different values of $-2<\b<-1,$ shown in Table \ref{mono}, clearly illustrates an example where $x_{2,n}$ decreases as $\b$ increases.
\end{itemize}
\end{remark}
\begin{table}[ht] \centering \caption{The values of $x_{1,5}$, $x_{2,5}$, $x_{3,5}$, $x_{4,5}$ and $x_{5,5}$ for $c=0.2$ and two different values of $\b$.}\label{mono}
\begin{tabular}{c|c|c|c|c|c}
$\b$&$x_{1,5}$&$x_{2,5}$&$x_{3,5}$&$x_{4,5}$&$x_{5,5}$\\
\hline
$-1.9$&$0.000000904$&$0.999651$&$2.006685$&$3.445917$&$6.1727379$\\
$-1.8$&$0.0000169239$&$0.999321$&$2.0144$ &$3.48815$&$6.24811$\\
\end{tabular}\end{table}

\section{Acknowledgments}
The authors are grateful to the anonymous referees for valuable comments. The research of KJ was supported by a Royal Society Newton Advanced Fellowship NAF$\backslash$R2$\backslash$180669. AJ and KJ thank their respective institutions for granting them research and development leave.


\end{document}